\documentclass[a4paper,12pt]{amsart}
\usepackage[english]{babel}
\usepackage[hyphens]{url} 
\usepackage{hyperref} 
\usepackage{amsmath} 
\usepackage{amssymb} 

\usepackage{amsthm}
\usepackage[all]{xy}
\usepackage{bussproofs}

\usepackage{color} 
\definecolor{boeuf}{rgb}{0.58,0,0.15}
\definecolor{brun}{rgb}{0.4,0.1,0.2}
\definecolor{bouteille}{rgb}{0.02,0.3,0.3}
\definecolor{marine}{rgb}{0.3,0.3,0.7}

\newcommand\ttt{\mathbf{t}}
\newcommand\eee{\mathbf{e}}

\newcommand\fl{\rightarrow}

\newcommand\lang{\mathcal{L}}

\newcommand\ee\epsilon 
\newcommand\imply\Rightarrow

\newdir{|>}{%
!/4.5pt/@{|}*:(1,-.2)@^{>}*:(1,+.2)@_{>}}

\mathcode`\:="603A     
\mathchardef\colon="303A  

\mathcode`\<="4268     
\mathcode`\>="5269     
\mathchardef\gt="313E  
\mathchardef\lt="313C  

\newtheorem{deff}{Definition}
\newtheorem{prop}[deff]{Proposition}

\newtheorem{Remark}[deff]{Remark}

\usepackage{gb4e}

\usepackage{mathptmx} 
\begin{document}
\title[Aristotle's square of opposition and Hilbert's quantifiers]{Aristotle's square of opposition in the light of\\ Hilbert's epsilon and tau quantifiers}

\author{Fabio Pasquali}
\address{\noindent Fabrio Pasquali\newline  CNRS IRIF\newline  B\^atiment Sophie Germain\newline  
8 place Aur{\'e}lie Nemours 
\newline 75013 Paris}
\email{pasquali@dima.unige.it}

\author{Christian Retor\'e} 
\address{Christian Retor\'e\newline  LIRMM universit\'e de Montpellier\newline  860 rue de Saint Priest b\^atiment 5\newline  34095 Montpellier cedex\newline\indent \tt \url{http://www.lirmm.fr/~retore} }
\email{christian.retore@lirmm.fr} 

\thanks{\it A preliminary work  has been presented by the second author at the 4th congress on the Square of opposition, Rome, 2014} 




\mbox{Work presented at}

\mbox{Aristotle 2400 Years World Congress, Thessaloniki, May 23--28 2016} 

\vspace*{3ex} 

\maketitle

\begin{abstract} 
Aristotle considered particular quantified sentences in his study of syllogisms and in his famous square of opposition. 
Of course, the logical formulas in Aristotle work were not modern formulas of mathematical logic, but ordinary sentences of natural language. 
Nowadays natural language sentences are turned into formulas of  predicate logic as defined by Frege, but, it is not clear that those Fregean sentences are faithful representations of natural language sentences. Indeed, the usual modelling of natural language quantifiers does not fully correspond to natural language syntax, as we shall see. 
This is the reason why Hilbert's epsilon and tau quantifiers (that go beyond usual quantifiers) 
have been used to model natural language quantifiers. 
Here we interpret Aristotle quantified sentences as formulas 
of Hilbert's epsilon and tau calculus. This yields to two potential squares of opposition and provided a natural condition holds, 
one of these two squares  is actually a  square of opposition i.e. satisfies the relations of contrary, contradictory, and subalternation. 
\end{abstract} 

\section{Aristotle I A E O sentences:\\ the standard interpretation and its inadequacies} 

When studying syllogisms and in his square of opposition, Aristotle considered 4 kinds of quantified sentences (see e.g. \cite{AristotleIEP}): 

\begin{tabular}{lll} 
A& Every S is P&	Universal Affirmative \\ 
I &Some S is P &	Particular Affirmative\\ 
E& No S is P&	Universal Negative \\ 
O&  Not all S are P & Particular negative \\   
\end{tabular} 

O formulas are different from the others: as observed e.g. in \cite{DelfittoVender2010ggag} no human language possesses a single word for saying "not every" or "not all" ---  at least in the language for which this question was asked. Nowadays,  O formulas are expressed by  \emph{Some S are not P}, and this induces a change of focus, from the global set of S entities (which gets out of P), to some particular entities in S (that are out of P):
\begin{exe} 
\ex Not all students may park their car on the campus. 
\ex Some students may not park their car on the campus. 
\end{exe} 

\subsection{Standard Montagovian modelling of quantifiers} 

Aristotle used to express quantified formulas in natural language, and, from a syntactic point of view,  there are two different kinds of quantifiers:

\begin{exe}
\ex  Something happened to me yesterday.  
\ex Some girls give me money.
\end{exe} 

The first kind has an \emph{explicit} domain, the restriction, which is a one-place predicate; the second kind of quantifiers does not have an explicit restriction, but   they all have an \emph{implicit} domain: because of their nature  (e.g. \emph{everything} does not range over human beings) or because of the linguistic and extra-linguistic context. As opposed to Frege's quantifiers for predicate logic natural language quantifiers do not range over a single sorted universe of  first order entities (for instance \emph{everything} may range over higher order objects, like $\ttt$ in : "He believes everything").

In the standard compositional, formal, and computable view of semantics  introduced by Montague \cite{Mon70,MootRetore2012LCG} formulas are represented 
by first order lambda terms over two base types $\eee$ (entities) and $\ttt$ (propositions) as Church did. \cite{Church1940types}
N-ary predicates are constants of the lambda calculus of type: $\eee\fl(\cdots\fl(\eee\fl\ttt)$ (N times $\eee$). 
Binary connectives $\&,\supset,\lor,...$ are constants of type $\ttt\fl(\ttt\fl\ttt)$. 
Quantifiers apply to predicates  to give a proposition, thus their type is 
$(\eee\fl\ttt)\fl\ttt$
and when there is a restriction to a class (e.g. \emph{some}) the quantifier applies to two predicates to give a proposition, thus its type is 
$(\eee\fl\ttt)\fl(\eee\fl\ttt)\fl\ttt$. 
Here are their precise definitions. 

\newcommand{\type}[1]{^{#1}}

\begin{center} 
\textit{something} \quad 
$ \exists\type{(e\fl t)\fl t}$ of type 
  $(e\fl t) \fl t)$

\textit{everything} \quad 
$ \forall\type{(e\fl t)\fl t}$ of type 
  $(e\fl t) \fl t)$

\textit{some} \quad 
 $\lambda P\type{e\fl t}\  \lambda Q\type{e\fl t}\  
(\exists\type{(e\fl t)\fl t}\  (\lambda x\type{e}  (\&\type{t\fl (t\fl t)} (P\ x) (Q\ x))))$ 
of type  $(e\fl t)\fl ((e\fl t) \fl t)$

\textit{every} \quad 
 $\lambda P\type{e\fl t}\  \lambda Q\type{e\fl t}\  
(\forall\type{(e\fl t)\fl t}\  (\lambda x\type{e}  (\supset\type{t\fl (t\fl t)} (P\ x) (Q\ x))))$ 
of type  $(e\fl t)\fl ((e\fl t) \fl t)$
\end{center}

\subsection{Inadequacies of the standard solution} 

There are basically three problems with this approach (see e.g. \cite{Retore2014fg}): 
\begin{enumerate} 
\item This view of quantifiers does not distinguish any exchange between the main predicate (rheme) and the domain predicate (theme). Alas it is clear that the two statements have a different meaning as their focuses are different. 
\begin{exe} 
\ex \label{crooks} 
\begin{xlist} 
\ex Some politicians are crooks.
\ex ?? Some crooks are politicians. 
\end{xlist}
\ex \label{employees} 
\begin{xlist} 
\ex Some students are employees. 
\ex Some employees are students. 
\end{xlist} 
\end{exe} 
\item The correspondence between syntax and semantics, that is the core of compositional semantics  is lost. Indeed, in the following example the underlined predicate does not correspond to a syntactic subtree (i.e. a constituent, a phrase), the syntactic structure and the semantic structure are different. 
\begin{exe}
\ex 
Keith composed some hits. 
\ex 
syntax (Keith  (composed  (some (hits))))
\ex 
semantics: (some (hits)) \underline{($\lambda x.$ Keith composed $x$)}
\end{exe} 
\item Interpretation of noun phrases and nominal phrases is difficult: indeed the interpretation of a noun phrase should be an entity, an individual, a term, i.e. should be of type $\eee$, and in case the main predicate is absent, as in nominal sentences, the standard approach  is not satisfactory, because it yields a type raised entity i.e. a term of type $(\eee\fl\ttt)\fl\ttt$. 
We indeed need a reference,  before the utterance of the main predicate (if any): 
\begin{exe}
\ex A goat! 
\ex  My thumb instead of an onion. (S.~Plath) 
\end{exe} 
\end{enumerate} 

\section{Hilbert's epsilon and tau quantifiers} 

Hilbert introduced in 
\cite{HilbertGrundlagen1922}
a logic with \emph{in situ} 
\footnote{as one says about interrogative pronouns in Chinese that do not move to the beginning of the sentence}  quantifiers in order to establish the consistency of arithmetics and analysis. 

A term (of type individual) $\epsilon_x F(x)$ is associated with every formula with  $F(x)$ and the operator $\epsilon$ binds all the free occurrences of $x$ in $F$ --- such operators are called subnectors in \cite{CF58}. The meaning is that as soon as an entity enjoys $F(\_)$ the term $\epsilon_x F(x)$ enjoys $F(\_)$: 
$$F(\epsilon_x F(x))\equiv \exists x.\ F(x)$$ 

Writing $x$ in $F(x)$ is unnecessary, indeed the variable $x$ may not appear in $F$ and their could be other variables. But it is a convenient notation, and furthermore $x$ usually appears in $F$.

There is a dual subnector: $\tau_x F(x)$ for which $$F(\tau_x F(x))\equiv \forall x.\ F(x)$$ 
Terms and formulae are defined by mutual recursion:
\begin{itemize} 
\item 
Any  constant in $\lang$ is a term.
\item 
Any  variable in $\lang$ is  a term.
\item 
$f(t_1,\ldots , t_p)$ is a term provided each $t_i$ is a term and $f$ is a function symbol of $\lang$ of arity $p$ 
\item 
$\epsilon_x A$ is a term if $A$ is a formula and $x$ a variable ---  any free occurrence of $x$ in $A$ is bound by $\epsilon_x$ 
\item 
$\tau_x A$ is a term if $A$ is a formula and $x$ a variable --- any free occurrence of $x$ in $A$ is bound by $\tau_x$ 
\item 
$s = t$ is a formula whenever $s$ and $t$ are terms. 
\item 
$R(t_1,\ldots,t_n)$  is a formula  provided each $t_i$ is a term and $R$ is a relation symbol of $\lang$ of arity $n$ 
\item 
$A \& B$, $A \lor B$, $A \Rightarrow  B$,  $\lnot A$ when $A$ and $B$ are formulae. 
\end{itemize} 

What are the proof rules for these termes, i.e. for quantification, in the epsilon calculus? 
They are the usual rules for quantification. 
From a  proof of $P(x)$ whose hypotheses do not concern $x$ (there are no free occurrences of $x$ in any hypothesis) infer  $P(\tau_x. P(x))$.

From a  proof of $P(x)$ whose hypotheses do not concern $x$ (there are no free occurrences of $x$ in any hypothesis) infer  $P(\tau_x. P(x))$.
From $P(\tau_x. P(x))$ one may infer $P(t)$ for any term $t$. 
From $P(t)$ infer $P(\epsilon_x P(x))$.  (the elimination rule for epsilon is slightly tricky we prefer not to give it, however the first rule together with the second one or the third one are enough) 

It is common (and sound) to add, in this variant of \emph{classical} logic, a rule $\epsilon x. F(x)=\tau x. \lnot F(x)$ . Indeed, the above deduction rules above shows that 
$P(\epsilon_x P(x)) \equiv \lnot \lnot P(\tau_x \lnot P(x))$ 
and 
$P(\tau_x P(x)) \equiv \lnot \lnot P(\epsilon_x \lnot P(x))$. 
 Hence one of the two subnectors $\epsilon$ and $\tau$ is enough, and usually people choose $\epsilon$. 
 
The quantifier free epsilon calculus is a strict conservative extension of first order logic. For more details on the epsilon calculus, quantification and natural language the reader is referred to \cite{HBvol2,epsilonIEP,espilonURL, APR2016,EgliHeusinger1995,Retore2014fg}. The reader  must be aware that there is no natural notion of model that interprets epsilon terms as usual first order terms or as choice functions: completeness fails, and this expected  first order and higher order logic cannot provide an account of the epsilon calculus.

%
%
%

\section{A I E O sentences with epsilon and tau} 

The I sentences,  \emph{some $S$ is  $P$} can be translated by 
$P(\epsilon x.\ S(x))$ as von Heusinger and other did \cite{EgliHeusinger1995,Retore2014fg}. 
It is 
not equivalent to any ordinary formula of first or higher-order predicate calculus, and in particular it is not equivalent to the standard: 
$\exists x.\ S(x)\&P(x)$. Nevertheless this formula has some relation with some usual formulas:
\footnote{ $A \vdash B$ means \emph{$A$ entails $B$} and $A \dashv\vdash B$ is a short hand for $A \vdash B$  and $B \vdash A$.}
$$\mathbf{P(\epsilon x.\ S(x))}\ \&\ \exists x. S(x) \vdash \exists x.\ (S(x)\&P(x))$$
$$\exists x. S(x)\ \&\ \forall y (S(y) \Rightarrow P(y)) \vdash \mathbf{P(\epsilon x. S(x))}$$ 

That way, all the three aforementioned  inadequacies vanish: 
\begin{enumerate} 
\item $S(\epsilon x. P(x)) \not\equiv  P(\epsilon x. S(x))$ hence the asymmetry is reflected in the semantic representation 
\item The interpretation is $((composed (\epsilon x. motet(x))) orlando)$ and follows the syntactic structure. 
\item $\epsilon x. goat(x)$ is a term of type $\eee$ as expected for a noun phrase --- the presupposition $goat(\epsilon x. goat(x)$ can be added in the context, because when someone says "a goat" there must be one, at least in his mind! 
\end{enumerate} 

The A formula  \emph{every $S$ is $P$} can be stated as the I formula \emph{some $S$ is $P$}  using $\tau$ instead of $\epsilon$: $P(\tau x.\ S(x))$
The two other formulas E and O are simply the negations of I and A respectively:

\noindent 
\begin{tabular}{lllll} 
A& Every S is P 
&  $P(\tau x. S(x))$&  $P(\tau S)$ & $A_{S,P}$ \\ 
I &Some S is P 
& $P(\epsilon x. S(x))$& $P(\epsilon S)$& $I_{S,P}$ \\ 
E& No S is P 
& $\lnot P(\epsilon x. S(x))\equiv 
\lnot P(\tau x. \lnot S(x))$ &$\lnot P(\epsilon S)$ & $E_{S,P}$ \\ 
O&  Some S are not P 
& $\lnot P(\tau x. S(x))\equiv 
\lnot P(\epsilon x. \lnot S(x))$ & 
$\lnot P(\tau S)$ 
& $O_{S,P}$\\   
\end{tabular} 

The two last columns define some useful abbreviations (because the predicates have no more than one free variable we do not need to specify the name of the variables that are bound by $\tau$ and $\epsilon$).  

\section{Two Hilbertian squares,\\ one of which is an Aristotelian square of opposition}

As an ongoing series of conference shows, the square of Aristotle is of a lasting interest. 
Since Aristotle use to study logic with natural language sentences and not with logical formulae
it is important that our new interpretation of quantified natural language sentences A E I O define a square of opposition for instance that contradictory sentences are logically contradictory. This possibility is mentioned in \cite{Slater2012square}. 

%
\begin{deff}\label{def} A square of opposition (see e.g. \cite{Slater2012square,sep-square} is given by four well formed formulas $A$, $E$, $I$, $O$ such that:
\footnote{ $A \vdash B$ means \emph{$A$ entails $B$} and $A \dashv\vdash B$ is a short hand for $A \vdash B$  and $B \vdash A$.}
\begin{itemize}
\item[i)] $A\dashv\vdash \neg O$ and $E\dashv\vdash \neg I$
\item[ii)] It is never the case that $\top\vdash A$ and $\top \vdash E$
\item[iii)] It is never the case that $I\vdash \bot$ and $E \vdash \bot$
\item[iv)] $A\vdash I$ and $E\vdash O$
\end{itemize}
\end{deff}
\begin{Remark}Note that condition $iii$ is redundant, since it is derivable from $i$ and $ii$. Also one of the conditions in $iv$ is deducible from the other and from $i$. 
\end{Remark}
Formulas $A$ and $O$ and formulas $E$ and $I$ are said to be contradictory, since one holds exactly when the other doesn't. By $(ii)$ formulas $A$ and $E$ can not be simultaneously true, thus they are said to be contraries. Conversely $I$ and $O$ cannot be simultaneously false, then they are called subcontraries. A proposition $X$ is said to be a subaltern of another proposition $Y$ when the truth of $X$ implies the truth of $Y$
(and some authors add that the falsity of $Y$ implies the falsity of $X$). 

The usual diagrammatical representation is
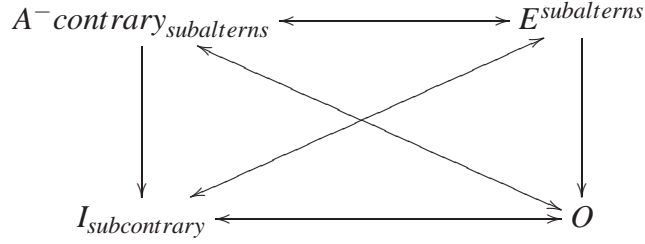
\[
\xymatrix{
A\ar@{<->}[ddrrr]\ar@{<->}[rrr]^-{contrary}\ar[dd]_{subalterns}&&&E\ar[dd]^{subalterns}\\
&&&\\
I\ar@{<->}[uurrr]\ar@{<->}[rrr]_{subcontrary}&&&O
}
\]


The motivating idea is that a term of the form $\ee_P$ is a witness of the fact that $P$ holds, so it makes sense to think at $Q(\ee_P)$ as "Some $Q$ are $P$". Nevertheless, since the Hilbert's $\ee$-calculus is classic, i.e the law of excluded middle holds, we have that $\top\vdash P \lor \neg P$. Then we might expect that one term among $\ee_P$ and $\ee_{\neg P}$ is a witness "better" then the other. This is made more precise by the following proposition. 


We say that \emph{$P$ is bivalent with respect to $S$} whenever $P(\epsilon S)\vdash  P(\tau S)$ 
or $P(\tau S) \vdash P(\epsilon S)$ --- observe that in the square case $S$ is a unary predicates, so the definition of bivalence only involves sentences i.e. formulae without free variables. The term “bivalent” comes from the fact that when $P$ is bivalent 
(i.e. forall $x$ $\vdash P(x)$ or  $\vdash \lnot P(x)$) 
$P$ is bivalent wrt any $S$.

Denote by $Sq (S, P)$ the square obtained with the following figures: $$A_{S,P},\quad  I_{S,P},\quad  E_{S,P},\quad  O_{S,P}.$$

Denote by $\mathcal{S}(S,P)$ the square obtained with the following figures: $A_{S,P}$, $I_{S,P}$, $E_{S,P}$ and $O_{S,P}$. 

\begin{prop}
If $P$ is bivalent with respect to $S$ then 
 either $\mathcal{S}(S,P)$ or $\mathcal{S}(S,\neg P)$ is a square of opposition.
\end{prop}
That is to say that, using standard diagrams, in the Hilbert's $\ee$-calculus one of the following is a square of opposition
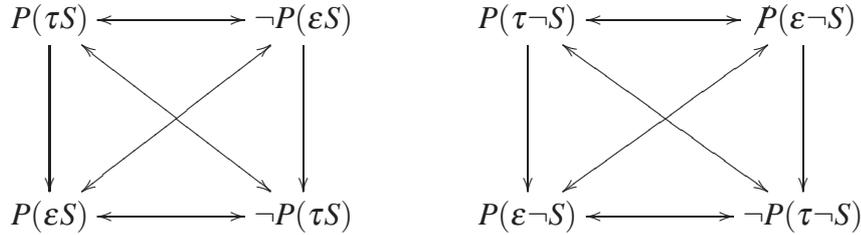
\[\xymatrix{
P(\tau S)\ar@{<->}[ddrr]\ar@{<->}[rr]\ar[dd]&&\neg P(\epsilon S)\ar[dd]\\
&&&\\
P(\epsilon S)\ar@{<->}[uurr]\ar@{<->}[rr]&&\lnot P(\tau S)
}\ \ \ 
\xymatrix{
P (\tau \lnot S)\ar@{<->}[ddrr]\ar@{<->}[rr]\ar[dd]&&\not P(\epsilon \lnot S)\ar[dd]\\
&&&\\
P(\epsilon \lnot S)\ar@{<->}[uurr]\ar@{<->}[rr]&&\neg P(\tau \lnot S)
}
\]
\begin{proof} Suppose $P(\tau S)\vdash P(\ee S)$.
\begin{description} 
\item[Contradictories] Diagonals are contradictories by construction.
\item[Subalterns] these coincide exactly with the assumption that $P (\tau S)\vdash P(\epsilon S)$.
\item[Contrary] We have to prove that $P(\tau S)$ and $\neg P(\epsilon S)$ cannot be both true. Suppose they are both true, i.e. 
\begin{enumerate}
\item \label{PtS}
$\top\vdash P(\tau S)$
\item  \label{PeS} 
 $\top\vdash\neg P(\epsilon S)$
 \end{enumerate} 
 By (\ref{PtS}) we have
$\top\vdash P(\tau S) \vdash P(\epsilon S)$ which contradicts (\ref{PeS}). 
\item[Subcontrary] we have to prove that $P(\epsilon S)$ and $\neg P(\tau S)$ cannot be both false. This comes for free by negating contraries.
\end{description} 

\noindent The other case, i.e id $P(\tau \neg S)\vdash P(\ee \neg S)$ is analogous.

\end{proof}

\section{Conclusion and future work} 

This study of the reformulation of the A E I O sentences with epsilon is particularly appealing both from a linguistic view point and from a logical one. 
Even for those very simple sentences, it is hard to grasp the intuitive meaning of the epsilon formulae that are not equivalent to first order formulae
and their relation to standard formulae. 

Regarding the square itself and its symmetries,  one may wonder how the square generalises with a typed epsilon that either acts directly on types or on predicates that apply to a given type. 

The proof that our square are square of opposition  relies on a strong version of \emph{tertium non datur}, and we wonder what happens with the intuitionistic variants of epsilon that were studied by Bell \cite{Bell1993mlq}.

\bibliographystyle{plainurl} 

\bibliography{bigbiblio}

\end{document}